\documentclass[11pt]{article}

\usepackage{a4wide}
\usepackage{amsmath,amssymb,amsfonts,amsthm, MnSymbol}
\usepackage{moreverb,rotating,graphics}
\usepackage[hidelinks]{hyperref}
\usepackage[T1]{fontenc}
\usepackage[utf8]{inputenc}
\usepackage{epsfig}
\usepackage[francais,english]{babel} 
\usepackage{framed,fancybox}
\usepackage{tikz,pgfplots}
\usepackage{graphicx}
\usepackage{caption}
\usepackage{subcaption}
\usepackage{mathrsfs}
\usepackage{esint}
\usepackage{dsfont}
\usepackage{enumerate}
\usepackage{import}
\usepackage[normalem]{ulem}
\usepackage{color}
\usepackage{afterpage}

\numberwithin{equation}{section} 

\newtheorem{theorem}{Theorem}[section]

\newtheorem{lemma}[theorem]{Lemma}

\newcommand\1{{\ensuremath {\mathds 1} }} 

\newcommand{\nm}{\noalign{\smallskip}}

\newcommand{\ds}{\displaystyle}

\def\C{{\mathbb C}}

\def\R{{\mathbb R}}

\def\br{{\bold r}}

\def\bs{{\bold s}}

\def\rd{{\mathrm{d}}}
\def\re{{\mathrm{e}}}
\def\ri{{\mathrm{i}}}

\def\Re{{\mathrm{Re }\, }}
\def\Im{{\mathrm{Im }\, }}

\def\cA{{\mathcal A}}

\def\cK{{\mathcal K}}
\def\cL{{\mathcal L}}
\def\cM{{\mathcal M}}

\def\cS{{\mathcal S}}

\def\rad{{\rm rad}}

\def\sqw{\hbox{\rlap{\leavevmode\raise.3ex\hbox{$\sqcap$}}$%
\sqcup$}}
\def\cqfd{\ifmmode\sqw\else{\ifhmode\unskip\fi\nobreak\hfil
\penalty50\hskip1em\null\nobreak\hfil\sqw
\parfillskip=0pt\finalhyphendemerits=0\endgraf}\fi}

\hypersetup{
    colorlinks,
    linkcolor={blue!50!blue},
    citecolor={red}
}
\renewcommand{\eqref}[1]{(\ref{#1})}

\newcommand\beq{\begin{equation}}
\newcommand\eeq{\end{equation}}
\newcommand\f[2]{\frac{#1}{#2}}
\newcommand\p{\partial}

\begin{document}

\title{Sub-wavelength focusing of acoustic waves in bubbly media\thanks{\footnotesize  
Hyundae Lee was supported by NRF-2015R1D1A1A01059357 grant.  Hai Zhang  supported by HK RGC grant ECS 26301016 and startup fund R9355 from HKUST.}}
\date{}

\author{
Habib Ammari\thanks{\footnotesize Department of Mathematics, 
ETH Z\"urich, 
R\"amistrasse 101, CH-8092 Z\"urich, Switzerland (habib.ammari@math.ethz.ch, brian.fitzpatrick@sam.math.ethz.ch).} \and Brian Fitzpatrick\footnotemark[2] \and David Gontier\thanks{\footnotesize CEREMADE, 
Universit\'e Paris-Dauphine, 
Place du Mar\'echal de Lattre de Tassigny, 
75775 Paris Cedex 16, France (gontier@ceremade.dauphine.fr).}
\and Hyundae Lee\thanks{\footnotesize  Department of Mathematics, Inha University,  253 Yonghyun-dong Nam-gu,  Incheon 402-751,  Korea (hdlee@inha.ac.kr).}  \and Hai Zhang\thanks{\footnotesize 
Department of Mathematics, 
 HKUST,  Clear Water Bay, Kowloon, Hong Kong (haizhang@ust.hk).}}

\maketitle

\begin{abstract}
The purpose of this paper is to investigate acoustic wave scattering by a large number of bubbles in a liquid at frequencies near the Minnaert resonance frequency. This bubbly media has been exploited in practice to obtain super-focusing of acoustic waves. Using layer potential techniques we derive the scattering function for a single spherical bubble excited by an incident wave in the low frequency regime. We then derive the point scatter approximation for multiple scattering by $N$ bubbles. We describe several numerical experiments based on the point scatterer approximation that demonstrate the possibility of achieving super-focusing using bubbly media.

\end{abstract}


\noindent {\footnotesize Mathematics Subject Classification
(MSC2000): 35R30, 35C20.}



\section{Introduction}
In this paper we are concerned with acoustic wave propagation in bubbly media. In particular, we are interested in the effect on the imaginary part of the Green's function of exciting  bubbles at frequencies near the Minnaert resonance frequency. It is well known that the focal spot size which limits the resolution limit in imaging is dependent on the imaginary part of the Green's function due to the Helmholtz-Kirchoff Theorem~\cite{ammari2013mathematical}. Physical experiments have shown that  it is possible to focus waves at the sub-wavelength scale by exploiting the strong scattering property of acoustic bubbles at their Minnaert resonance~\cite{Lanoy2015}.

\medskip

The Minnaert resonance is a quasi-static resonance in which the bubble size is much smaller than the wavelength of the incident wave. We use layer potential techniques to explicitly determine the acoustic surface potentials for a single bubble in the quasi-static regime, in which sound propagation can be neglected and the pressure field on the surface of the bubble may be considered a constant~\cite{Ammari2009_book}. Knowledge of the acoustic potential corresponding to the exterior domain enables us to the derive the scattering function for a single bubble, which has a similar expression to the scattering function found in the physical literature. The Minnaert resonance corresponds to the value at which the scattering function is maximized.

\medskip
In \cite{H3a}, a rigorous mathematical justification of the Minnaert resonance in the case of a single bubble in a homogeneous medium is established. The acoustic properties of the bubble are analyzed, and the approximate formula for the Minnaert resonance of an arbitrary shaped bubble is derived.  In \cite{phononic} the opening of a sub-wavelength phononic bandgap is demonstrated by considering a periodic arrangement of bubbles and exploiting their Minnaert resonance.  As shown in \cite{Ammari_Hai}, around  the Minnaert resonant frequency, an effective medium theory can be derived in the dilute regime. Furthermore, {above} the Minnaert resonant frequency, the real part of the effective {{modulus}} is negative and consequently, the bubbly fluid behaves as a diffusive media for the acoustic waves. Meanwhile, below the Minnaert resonant frequency, with an appropriate bubble volume fraction, a high contrast effective medium can be obtained, making the sub-wavelength focusing or super-focusing of waves achievable \cite{Ammari2015_a}.  {These properties show that the bubbly fluid functions like an acoustic metamaterial and indicate that a sub-wavelength bandgap opening occurs at the Minneaert resonant frequency \cite{Leroy2009} .} We remark that such behavior is rather analogous to the coupling of electromagnetic waves with plasmonic nanoparticles, which results in effective negative or high contrast dielectric constants for frequencies near the plasmonic resonance frequencies \cite{deng,matias1,matias2}.

\medskip
By generalizing the approach of \cite{H3a} to the case of $N$ bubbles all having the same radius, we obtain an $N \times N$ block matrix in which the diagonal terms represent bubble self-interaction, while the off-diagonal terms represent the interaction between different bubbles. Letting the size of the bubbles go to zero we derive the point scatterer approximation which reduces the problem to an $N \times N$ linear system, which is recognizable as the classical Foldy-Lax formulation for a scattered field due to a collection of point scatterers~\cite{Devaud2008}.

\medskip
We numerically simulate the system of point scatterers to analyze the super-focusing phenomenon observed in~\cite{Lanoy2015} and to simulate a time reversal experiment that demonstrates that it is possible to spatially localize a time reversed signal to a very high degree of accuracy in a bubbly media.

\medskip

Due to the Helmholtz-Kirchhoff Theorem, the sharper the imaginary part of the Green's function, the smaller the focal size and the higher the resolution achieved \cite{ammari2013mathematical, Ammari2015, Ammari2015_a}. To better understand the effect of the bubbly media on the imaginary part of the Greens function, we excite the point scatterers over a wide range of frequencies. By averaging the imaginary part of the Green's function function over a range of frequencies, we demonstrate that sub-wavelength focusing is achievable in the region slightly below the Minnaert resonance frequency.

\medskip


\medskip

The paper is organized as follows. In Section~\ref{sec_single_bubble} we formulate the scattering problem for a single bubble and derive the scattering function. In Section~\ref{sec_sys_of_N_bubbles} we derive the point scatterer approximation for a system of $N$ bubbles. In Section~\ref{sec_super_focusing} we demonstrate  the super-focusing phenomenon by numerically showing that the imaginary part of the Green's function becomes sharper when the incident wave frequency is close to the Minnaert resonance frequency. 


\section{Single bubble} \label{sec_single_bubble}

The purpose of this section is to justify the usual scattering function that one can find in the literature for bubbles. We explain the physical argument given in~\cite{Devaud2008} using a layer potential approach~\cite{Ammari2009_book}. In what follows we will only consider spherical bubbles. This is not a strong assumption for the phenomena we are interested in, and it already contains all the properties of a more complicated model; see \cite{H3a}. We first consider a single bubble, which will be represented by a sphere $D$ of radius $R$, centered at the origin:
\[
	D := \left\{ x \in \R^3, \ | x | \le R \right\}.
\]

We denote by $\rho_b$ and $\kappa_b$ the density and the bulk modulus of the air inside the bubble, respectively, and let $\rho_w$ and $\kappa_w$ be the corresponding parameters for the background media $\R^3 \backslash D$. The scattering problem can be modeled by the following equations:
\beq \label{eq-scattering}
\displaystyle
\left\{
\begin{array} {ll}
&\ds \nabla \cdot \f{1}{\rho_w} \nabla  u+ \frac{\omega^2}{\kappa_w} u  = 0 \quad \mbox{in } \R^3 \backslash D, \\
\nm
&\ds \nabla \cdot \f{1}{\rho_b} \nabla  u+ \frac{\omega^2}{\kappa_b} u  = 0 \quad \mbox{in } D, \\
\nm
&\ds u_{+} -u_{-}  =0   \quad \mbox{on } \partial D, \\
\nm
&  \ds \f{1}{\rho_w} \f{\p u}{\p \nu} \bigg|_{+} - \f{1}{\rho_b} \f{\p u}{\p \nu} \bigg|_{-} =0 \quad \mbox{on } \partial D,\\
\nm
&  u^s:= u- u^{i}  \,\,\,  \mbox{satisfies the Sommerfeld radiation condition.}
  \end{array}
 \right.
\eeq
We introduce four auxiliary parameters to facilitate our analysis:
$$
c_w = \sqrt{\frac{\kappa_w}{\rho_w}}, \, \, c_b = \sqrt{\frac{\kappa_b}{\rho_b}}, \, \, k_w =\f{\omega}{c_w}, \,\, k_b = \f{\omega}{c_b}.
$$
$c_w$ and $c_b$ represent the speed of sound in the background media and in the bubble, respectively. We also introduce the dimensionless contrast parameter
$$
\delta = \f{\rho_b}{\rho_w},
$$
and note that for water and air in regular conditions we have that $\delta \approx 1.2 \times 10^{-3}$.
\subsection{Layer potentials}

We use layer potentials to represent the solution to the scattering problem (\ref{eq-scattering}). Let the single layer potential $\mathcal{S}_{D}^{k}$ associated with the domain $D$ and the wavenumber $k$ be defined by
$$
\mathcal{S}_{D}^{k} [\psi](x) =  \int_{\p D} G(x, y, k) \psi(y) d\sigma(y),  \quad x \in  \p {D},
$$
where
\beq \label{eq:def:GF}
G(x, y, k)= G(x-y, k)=- \f{e^{ik|x-y|}}{4 \pi|x-y|},
\eeq
is the Green function for the Helmholtz equation in $\R^3$, subject to the Sommerfeld radiation condition.
We also define the boundary integral operator $\mathcal{K}_{D}^{k, *}$ by
$$
\mathcal{K}_{D}^{k, *} [\psi](x)  = \int_{\p D } \f{\p G(x, y, k)}{ \p \nu(x)} \psi(y) d\sigma(y) ,  \quad x \in \p D. 
$$

Then the solution $u$ can be written as 
\beq \label{Helm-solution}
u(x) = \left\{
\begin{array}{lr}
u^{in} + \mathcal{S}_{D}^{k_w} [\psi_w], & \quad x \in \R^3 \backslash \bar{D}, \\
\mathcal{S}_{D}^{k_b} [\psi_b] ,  & \quad x \in {D},
\end{array}\right.
\eeq
for some surface potentials $\psi_w, \psi_b \in  L^2(\p D)$. 
Using the jump relations for the single layer potentials, it is easy to derive that $\psi_w$ and $\psi_b$ satisfy the following system of boundary integral equations:
\beq \label{eq-boundary}
\mathcal{A}_D(\omega, \delta)[\Psi] =F,  
\eeq
where
\[
\mathcal{A}_D(\omega, \delta) = 
 \begin{pmatrix}
  \mathcal{S}_D^{k_b} &  -\mathcal{S}_D^{k_w}  \\
  \f{1}{\delta}(-\f{1}{2}Id+ \mathcal{K}_D^{k_b, *})& -( \f{1}{2}Id+ \mathcal{K}_D^{k_w, *})
\end{pmatrix}, 
\,\, \Psi= 
\begin{pmatrix}
\psi_b\\
\psi_w
\end{pmatrix}, 
\,\,F= 
\begin{pmatrix}
u^{in}\\
\f{\partial u^{in}}{\partial \nu}
\end{pmatrix}.
\]

One can show that the scattering problem (\ref{eq-scattering}) is equivalent to the boundary integral equations (\ref{eq-boundary}).


\subsection{The breathing approximation}

\subsubsection{The matrix $A_D$}
Note that on the surface of the bubble, we have $u^{in}(x)=u^{in}(0) + O(k_wR)$, and $\f{\partial u^{in}}{\partial \nu} = O(k_wR)$. 
If the wavelength is much larger than the size of the bubble, \textit{i.e.}, $k_w R \ll 1$, then we may approximate $u^{in}$ by $u^{in}(0)$ and $\f{\partial u^{in}}{\partial \nu}$ by $0$ on the surface of the bubble, respectively. In this case, we may look for solutions to (\ref{eq-boundary}) with the right-hand side given by
$
\begin{pmatrix} \1_{\partial D} \\ 0 \end{pmatrix}. 
$
This is the so-called \textit{breathing approximation}~\cite{Devaud2008}. 

\medskip

Thanks to the spherical symmetry of $D$, the operator $\mathcal{A}_D$
can be diagonalized by spherical harmonic functions. Denote by $X = L^2(\p D)\times L^2(\p D)$.
We define the finite-dimensional subspace $X_1$ of $X$ by
\[
	X_1 := {\rm Vect} \left\{ \begin{pmatrix} \1_{\partial D} \\ 0 \end{pmatrix}, \  \begin{pmatrix} 0 \\ \1_{\partial D}  \end{pmatrix} \right\}.
\]
The vector space $X_1$ is of dimension $2$. Using the breathing approximation, we only need to solve the following equation
$$
\mathcal{A}_D [\Psi] = \begin{pmatrix} \1_{\partial D} \\ 0 \end{pmatrix},
$$
for $\Psi \in X_1$.

\medskip

Let $A_D$ be the $2 \times 2$ matrix representing the operator $\cA_D$ on $X_1$. The values of the components of $A_D$ are deduced from the next lemma (see Section~\ref{sec:proof:SD_KD[1]} for the proof).
\begin{lemma} \label{lem:SD_KD[1]}
	It holds that
	\[
		\cS_D^k[ \1_{\partial D} ] = - \dfrac{\re^{\ri k R}}{k} \sin ( k R),
		\]
		and
		\[
		\left( \cK_D^k \right)^*[ \1_{\partial D} ] = \dfrac{\sin (kR)}{kR} \re^{\ri R k} - \frac{\re^{2 \ri kR}}{2}.
	\]
\end{lemma}

By introducing the dimensionless parameters $x_b = k_b R$ and $x_w = k_w R$, we obtain
\[
	A_D = \begin{pmatrix}
		- \dfrac{\re^{\ri x_b}}{x_b} \sin ( x_b) & \dfrac{\re^{\ri x_w}}{x_w} \sin ( x_w)  \\
		& \\
		\dfrac{\re^{\ri x_b}}{\delta} \left( \dfrac{\sin (x_b )}{x_b } - \cos \left(x_b \right) \right) & - \re^{\ri x_w} \left( \dfrac{\sin (x_w)}{x_w} - \ri  \sin\left( x_w  \right) \right)
	\end{pmatrix}.
\]

Denote $\psi_{b/w} = C_{b/w} \1_{\partial D}$ where $(C_b, C_w) \in \C^2$. 
Then we have
\begin{equation} \label{eq:solve1bubble}
	A_D \begin{pmatrix} C_b \\ C_w \end{pmatrix} = \begin{pmatrix} \frac{u^{in}(0)}{R} \\ 0 \end{pmatrix}.
\end{equation}

The determinant of $A_D$ is
\[
	\det(A_D) = \re^{\ri (x_b + x_w)} \sin(x_w) \left( \dfrac{\sin(x_b)}{ x_b} \left( \dfrac{1}{x_w} - \ri \right)  - \dfrac{1}{\delta x_w} \left( \dfrac{\sin(x_b)}{x_b} - \cos(x_b) \right)\right),
\]
and the inverse of $A_D$ is
\[
	\left( A_D \right)^{-1} = \dfrac{-1}{\det(A_D)} 
	\begin{pmatrix}
		\re^{\ri x_w} \left( \dfrac{\sin (x_w)}{x_w} - \ri  \sin\left( x_w  \right) \right) 
		& \dfrac{\re^{\ri x_w}}{x_w} \sin ( x_w) 
		 \\
		& \\
		\dfrac{\re^{\ri x_b}}{\delta} \left( \dfrac{\sin (x_b )}{x_b } - \cos \left(x_b \right) \right) & \dfrac{\re^{\ri x_b}}{x_b} \sin ( x_b) 
	\end{pmatrix}.
\]
Together with~\eqref{eq:solve1bubble}, we obtain that (we assume that $\sin(x_b) \neq 0$ for simplicity)
\begin{equation} \label{eq:Cwexplicit}
	C_w = \left[ \dfrac{ \re^{-\ri x_w} \left( \dfrac{1}{x_b } - \cot \left(x_b \right) \right) }{\sin(x_w) \left( \dfrac{1}{x_w} \left( \dfrac{1}{x_b} - \cot(x_b) \right) - \dfrac{\delta}{x_b} \left( \dfrac{1}{x_w} - \ri \right) \right)} \right] \frac{u^{in}(0)}{R} .
\end{equation}



\subsubsection{Scattering function}

 The pressure radiated by the bubble is ${u_\rad} := \cS_D^{k_w}[\psi_w]$. Since $\psi_w = C_w \1_{\partial D}$, we obtain ${u_\rad} = C_w  \cS_D^{k_w}[\1_{\partial D}]$. The value of this last function is given by the next Lemma (see also~\eqref{eq:explicitS1}).
\begin{lemma} \label{lem:widetildeS1}
	It holds that
	\[
		\forall x \in \R^3 \setminus D, \quad \cS_D^{k_w}[\1_{\partial D}](x) = - \dfrac{R \sin ( k_w R)}{k_w} \dfrac{ \re^{\ri k_w | x |}}{ | x |}.
	\]
\end{lemma}
Together with~\eqref{eq:Cwexplicit}, we deduce that the value of ${u_\rad}$ on the boundary $\partial D$ is (we assume that $\sin(x_w) \neq 0$)
\[
	u_\rad(|x| = R) = - C_w R \left( \dfrac{ \sin ( x_w) }{x_w} \re^{\ri x_w} \right)
	= \left( \dfrac{ x_b \cot \left(x_b \right) - 1}{  1 - \delta - x_b \cot (x_b )   + \ri \delta \dfrac{ c_b}{c_w } x_b} \right) u^{in}(0).
\]
Typical orders of magnitude are $\delta \approx 10^{-3}$ and $\delta \dfrac{c_b}{ c_w} \approx 2 \times 10^{-4}$. The dimensionless function
\begin{equation} \label{eq:def:fs}
	f_s(x_b) = \dfrac{ R(1 - x_b \cot \left(x_b \right) ) }{  x_b \cot (x_b ) - 1 + \delta - \ri \delta \dfrac{c_b}{c_w } x_b},
\end{equation}
is called the \textit{scattering function}. It links the value of the incoming pressure on  the boundary of the bubble to the value of the radiated pressure on the same boundary. Using the scattering function, we have the following monopole approximation for the scattered field
\begin{equation} \label{monopole}
	u_\rad(x) =	f_s(x_b)u^{in}(0) \dfrac{ \re^{ i k|x|} }{|x|}.
\end{equation}

\medskip

Let us study the function $f_s$. The \textit{Minnaert resonance} corresponds to the point $x_b = x_M$ such that
\[
	1 - \delta - x_M \cot (x_M) = 0
	\quad \text{or} \quad
	\tan (x_M) = \dfrac{x_M}{1 - \delta}.
\]
This point corresponds to the frequency $x_M \approx 0.06$. At this particular frequency, it holds that
\[
	f_s(x_M) = \ri \dfrac{c_w}{c_b x_M} 
	\quad \text{with} \quad
	\dfrac{c_w}{c_b x_M}  \approx 70.
\]
Hence, at the Minnaert frequency, the pressure gains a phase of $\re^{\ri \pi/2}$, and the amplitude is multiplied by a factor of $70$. For instance, for a bubble of radius $R = 10^{-3} \ m$, we obtain a frequency $\omega_M/(2 \pi) = \dfrac{c_b x_M}{2 \pi R } \approx 3300 \ \text{Hz}$, which is a perfectly audible sound. Note that $x_M$ \textit{does not} correspond to any extremal point of $f_s$. This is however a vanishing point for the real part. We can also consider the frequency corresponding to $x = x_0$ such that
\[
	\left| f_s(x_0) \right| = \sup_{x \in [0,0.1]} \left| f_s(x) \right|.
\]
Actually, it holds that $\left| x_M - x_0 \right| \approx 6.10^{-4}$, so that $x_M$ is indeed a very good (and easily computable) approximation of $x_0$. 

\medskip

\begin{figure}[h!]
    \centering
     \includegraphics[scale=0.48]{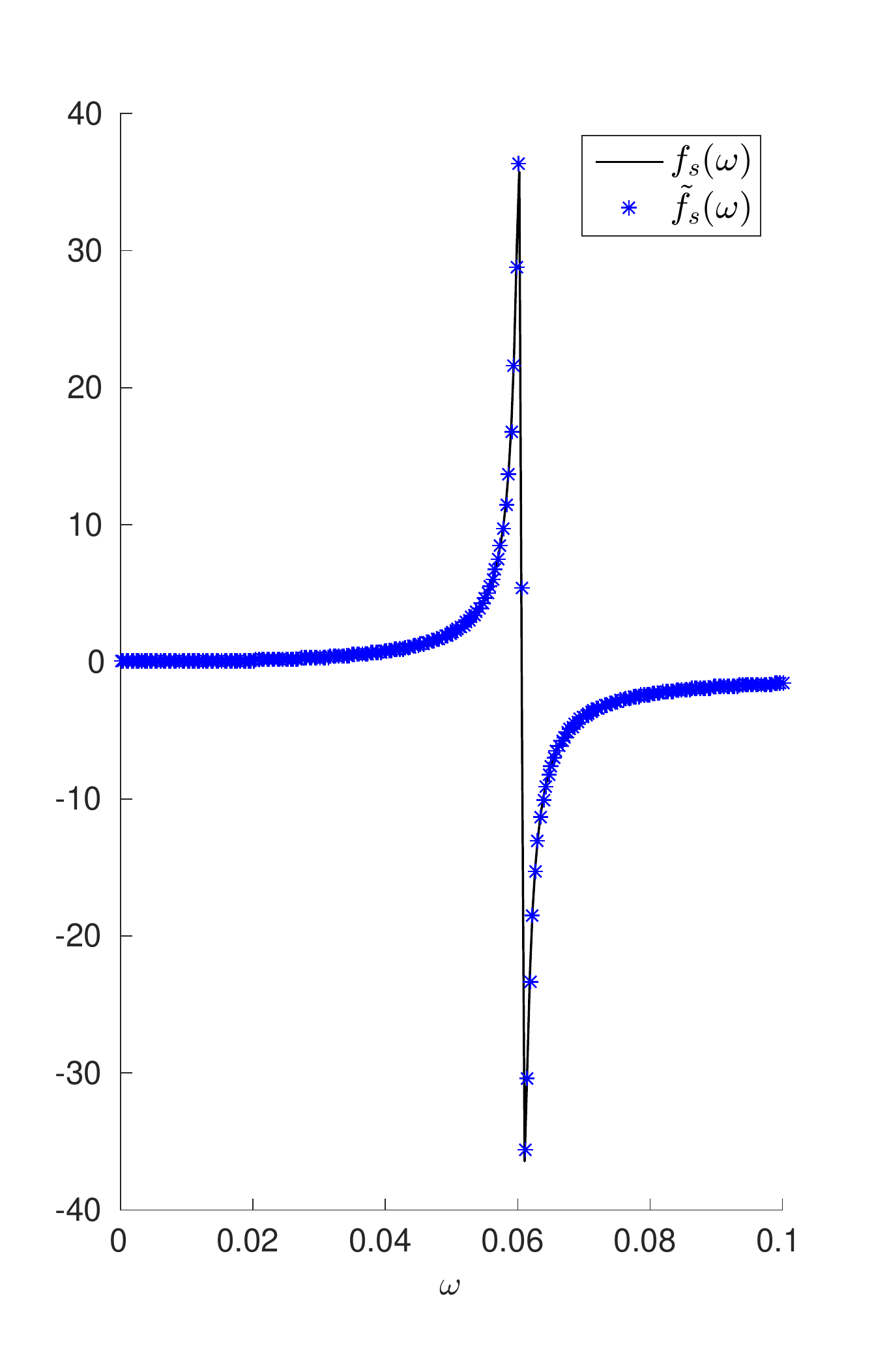}
     \includegraphics[scale=0.48]{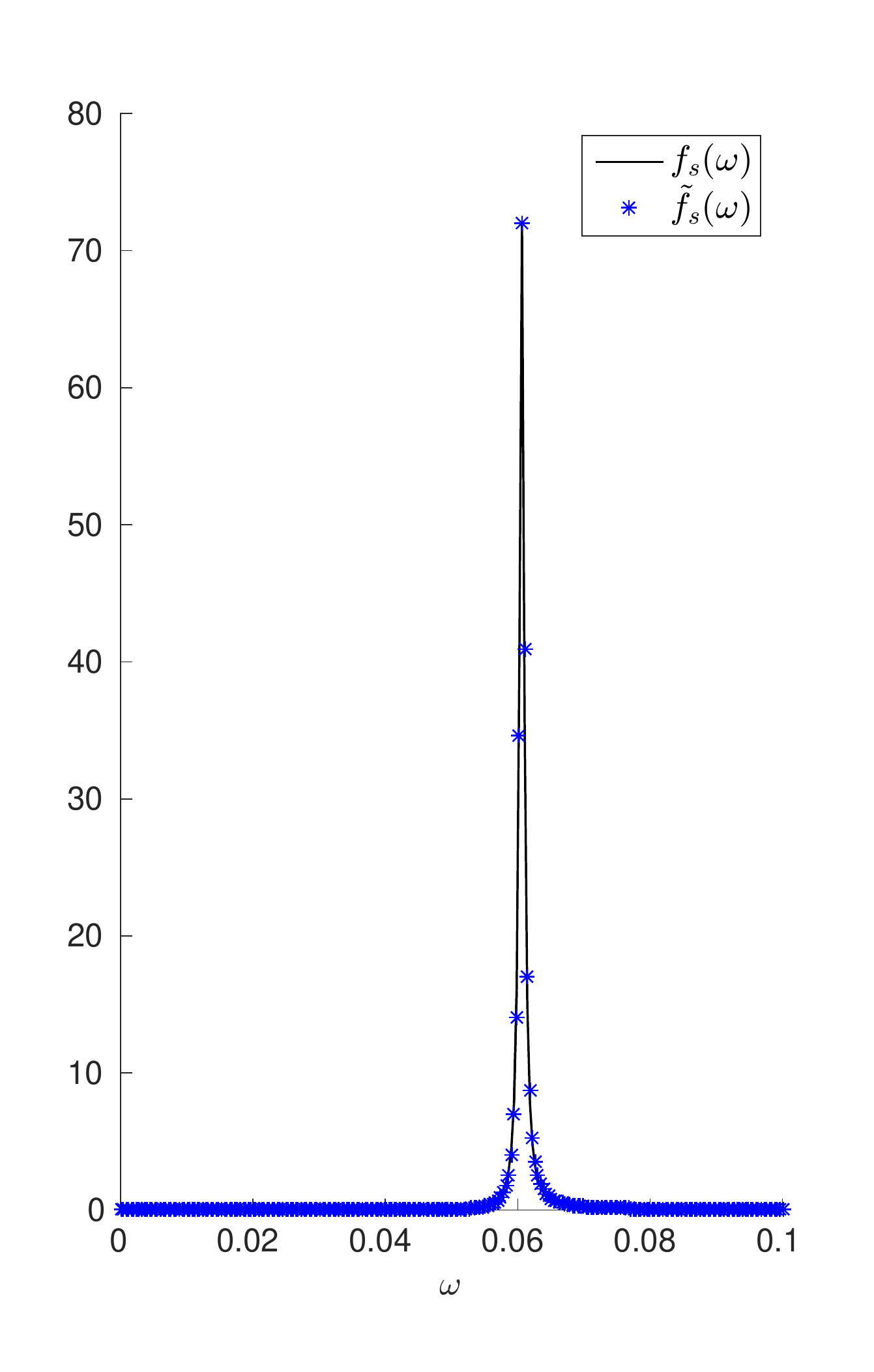}     
    \caption{The real part and the imaginary part of $f_s$ (solid black line) and $\widetilde{f_s}$ (blue crosses).}
    \label{fig:real_imag_part_fs}
\end{figure}

In the physical literature, the scattering function is usually approximated by the simpler function  
\begin{equation} \label{eq:fsTilde}
	f_s(x_b) \approx \widetilde{f_s}(x_b) := \dfrac{1}{\left( \frac{x_M}{x_b} \right)^2 - 1 - \ri \left( \frac{c_b}{c_w} \right)x_b}.
\end{equation}
The functions $f_s(\omega)$ and $\tilde{f}_s(\omega)$ are shown in Figure~\ref{fig:real_imag_part_fs}. Numerically, we find that $\sup_{x \in [0, 0.1]} \left| f_s - \widetilde{f_s} \right| \approx 0.47$, and hence $\widetilde{f_s}$ can be considered a highly accurate approximation of $f_s$.

\section{System of N bubbles} \label{sec_sys_of_N_bubbles}

We now consider a set of $N$ disjoint bubbles. In order to keep the notation simple, we will assume that all the bubbles have the same radius $R$. This is not a very restrictive assumption, and this case contains effects we are interested in. The volume occupied by bubble number $j$ is
\[
	D_j := \left\{ x \in \R^3, \ | x - x_j | \le R \right\}, \quad 1 \le j \le N.
\]
The distance between the centers of bubble $i$ and bubble $j$ is $d_{ij} = |x_i - x_j|$. The non-intersecting condition implies that $d_{ij} \ge 2R$ for all $i \neq j$. We finally denote by $D := \bigcup_{1 \le j \le N} D_j$ the volume occupied by all the bubbles.

\medskip

We want to solve the Helmholtz equation $(\Delta + k^2) u = 0$ subject to boundary conditions analogous to the ones in~\eqref{eq-scattering}. We make the ansatz
\begin{equation} \label{eq:ansatz2bubbles}
	u = \left\{
		\begin{array}{ll}
		 \left(\sum\limits_{j=1}^N \cS_{D_j}^{k_w}[\phi_{j,w}] \right)   + u^{in}, & x \in \R^3 \setminus \overline{D}, \\
		 & \\
		\cS_{D_j}^{k_b}[\phi_{j,b}] & x \in D_j, \quad 1 \le j \le N.\\
		\end{array} \right. 
\end{equation}
By introducing $X := \Pi_{j=1}^N L^2(\partial D_j) \times  L^2(\partial D_j)$ , this leads to a system of equations of the form
\begin{equation} \label{eq:Nbubbles}
	\cA_{D_1, \ldots D_N}
	\begin{pmatrix}
		\psi_{1,b} \\
		\psi_{1,w} \\
		\vdots \\
		\psi_{N,b} \\
		\psi_{N,w} \\
	\end{pmatrix}
	=
	\begin{pmatrix}
		u^{in} \Big|_{\partial D_1}  \\
		\f{\p u^{in}}{\p \nu_1}\Big|_{\partial D_1}  \\
		\vdots \\
		u^{in} \Big|_{\partial D_N} \\
		\f{\p u^{in}}{\p \nu_N}\Big|_{\partial D_N}  \\
	\end{pmatrix},
\end{equation}
where $\cA_{D_1, \ldots, D_N}$ is acting on $X$ and where $\Phi := (\psi_{1,b}, \psi_{1, w}, \ldots, \psi_{1,b}, \psi_{1, w})^T \in X$ with $T$ being the transpose. The operator $\cA_{D_1, \ldots, D_N}$ has the block diagonal form
\[
\cA_{D_1, \ldots, D_N} = \begin{pmatrix}
	\cM_1 & \cL_{1,2} & \cL_{1,3} & \ldots \\
	\cL_{2,1} & \cM_2 & \cL_{2,3} & \ldots \\
	\vdots & \vdots & \vdots & \vdots \\
	\cL_{N,1} & \cL_{N,2} & \ldots & \cM_N
\end{pmatrix},
\]
where (see also~\eqref{eq-boundary})
\[
	\cM_j :=  \begin{pmatrix}
		 \cS_{D_j}^{k_b} & - \cS_{D_j}^{k_w}  \\
		\dfrac{1}{\delta} \left( - \frac12 + \left(\cK_{D_j}^{k_b}\right)^* \right) & - \left( \frac12 + \left(\cK_{D_j}^{k_w}\right)^* \right) 
	\end{pmatrix},
\]
and, for $i \neq j$,
\[
	\cL_{i,j} = \begin{pmatrix}
		0 & - \cS_{D_i, D_j}^{k_w} \\
		 0 & - \cL_{D_i, D_j}^{k_w}
	\end{pmatrix}.
\]
Here, we introduced the operators $\cS_{D_i, D_j}^{k} : L^2(\partial D_j) \to L^2(\partial D_i)$ and $\cL^k_{D_i, D_j} : L^2(\partial D_j) \to L^2(\partial D_i)$, which describe the effect of bubble $j$ on bubble $i$. These operators are respectively defined by
\begin{equation*}
	\forall \phi \in L^2(\partial D_j) , \quad \cS_{D_i, D_j}^{k}\left[ \phi \right]  = \cS_{D_j}^k[\phi] \Big|_{\partial D_i}, 
	\quad \text{and} \quad
	\cL_{D_i, D_j}^{k}\left[ \phi \right]  = \f{\p}{\p \nu_i} \cS_{D_j}^k[\phi] \Big|_{\partial D_i}.
\end{equation*}

\medskip

 Such a system is usually called a \textit{system of $N$ scatterers} for the Helmholtz equation, and has been the starting point for a large number of articles.

\subsection{The point scatterer approximation for N-bubbles}
We derive the point scatterer approximation for the N-bubble system. For this, we assume the breathing approximation for each of the bubbles, \textit{i.e.}, the pressure field applied to each bubble is constant. More precisely, for each $1\leq i \leq N$, we define 
$u^{in}_i$  to be the total field incident on the $i$-th bubble $D_i$, and $u^s_i$ to be the field scattered from $D_i$. It is clear that
\begin{equation} \label{eq-11}
u^{in}_i(x) = u^{in}(x) + \sum_{j\neq i}u_j^s(x).
\end{equation}
Due to our assumption, $u^{in}_i$ can be approximated by a constant function across the interface $\partial D_i$. This is a good approximation when the bubbles are well-separated, \textit{i.e.}, their distance is much greater than their size. Under this approximation, we can derive a good approximate solution for the system (\ref{eq:Nbubbles}).
Indeed, the constant $u^{in}_i(x_i)$ can be used to approximate the pressure impinged on $D_i$. As a result, using (\ref{monopole}), the scattered field $u^s_i$ can be approximated by
\[
u^s_i(x) = -4\pi f_s(x_b) u^{in}_i(x_i)G(x, x_i, k_w)= -4\pi f_s(x_b) u^{in}_i(x_i)G(x-x_i, k_w),
\]
By taking $x=x_j$ in the above equation and using 
(\ref{eq-11}), we obtain the following systems of linear equations
for $u^{in}_i(x_i)$: 
\[
u^{in}_i(x_i) = u^{in}(x_i) + \sum_{j\neq i} -4\pi f_s(x_b) u^{in}_j(x_j)G(x_i-x_j, k_w), \quad 1\leq i \leq N.
\]
We define a matrix $M$ by setting
\begin{equation} \label{matrix-m}
M_{ij} =
\begin{cases}
1, & i = j, \\
4 \pi f^s G(x_i -x_j, k_w), & i \neq j.
\end{cases}
\end{equation}
Then the total field for the scattering problem by N bubbles can be approximated by
\[
u(x) = u^{in}(x)+ \sum_{1\leq i \leq N} -4\pi f_s(x_b) G(x -x_j, k_w)\sum_{1\leq j \leq N} (M^{-1})_{ij}u^{in}(x_j),
\]
where $(M^{-1})_{ij}$ denotes the $(ij)$-th component of the inverse  matrix $M^{-1}$. This is called the point interaction approximation. 
See also \cite{Ammari_Hai} for a rigorous justification of this approximation. 
\section{Super focusing} \label{sec_super_focusing}

We now perform a numerical investigation of the results that have been observed in~\cite{Lanoy2015} in order to assist us in developing a mathematical theory that explains the super-focusing phenomenon.



We take the density of water and air to be $\rho_w = 1 \times10^3 \ \text{kg m}^{-3}$, and $\rho_b = 1.2 \ \text{kg m}^{-3}$, respectively. We take the bulk modulus of water and air to be $\kappa_b = 2.07 \times10^9 \ \text{N m}^{-2}$, and $\kappa_b = 127 \times10^3 \ \text{N m}^{-2}$, respectively. This results in a speed of sound in water of $c_w \approx 1440 \ \text{m s}^{-2}$.

We set the common radius of the bubbles to be $R = 50 \ \mu \ \text{m}$. In light of these parameters, the Minnaert resonant frequency is $\omega_M \approx 59\times 2 \pi \ \text{kHz} \approx 390 \ \text{kHz}$. We are interested in exciting the bubbles with a reference frequency $\omega_R$ which is slightly below the Minnaert resonance frequency. Hence, we choose $\omega_R = 57.5 \times 2 \pi \ \text{kHz}$. We denote by by $\lambda_w = 2\pi/k_w$ the wavelength in water. Typical values for the wavelength are $\lambda_w \approx 0.025 \ m$.  The scattering function $f_s$ is taken to be
\[
	f_s(\omega) := \dfrac{1}{\left( \frac{\omega_M}{\omega} \right)^2 - 1 - \ri \epsilon}, \quad \text{with} \quad \epsilon = Rk_w.
\]
This expression is slightly different than~\eqref{eq:fsTilde}, although it behaves in a similar manner. It is worth noting that the parameters we have chosen essentially mirror those used in the experiments described in~\cite{Lanoy2015}.

\medskip

Our experimental setup is also largely in agreement with~\cite{Lanoy2015}. We consider a source, located at the origin $\bs = (0,0,0)^T$, that is surrounded by a cube containing uniformly distributed bubbles. The length of the cube is $L = 0.01 \ \text{m} \approx \lambda_w/2$. The gas volume fraction in the cube is $\Phi = 2\times 10^{-4}$, which corresponds to a collection of approximately $380$ bubbles. We assume that the water outside this cube does not contain any bubbles. Finally, four receivers (time-reversal mirrors) are located at $\br_1 = (0.02, 0, 0)^T$, $\br_2 = (-0.02, 0, 0)^T$, $\br_3 = (0, 0.02, 0)^T$ and $\br_4 = (0, -0.02, 0)^T$. 

\medskip

\begin{figure}[h!]
    \hspace{-2cm} \includegraphics[scale=0.5]{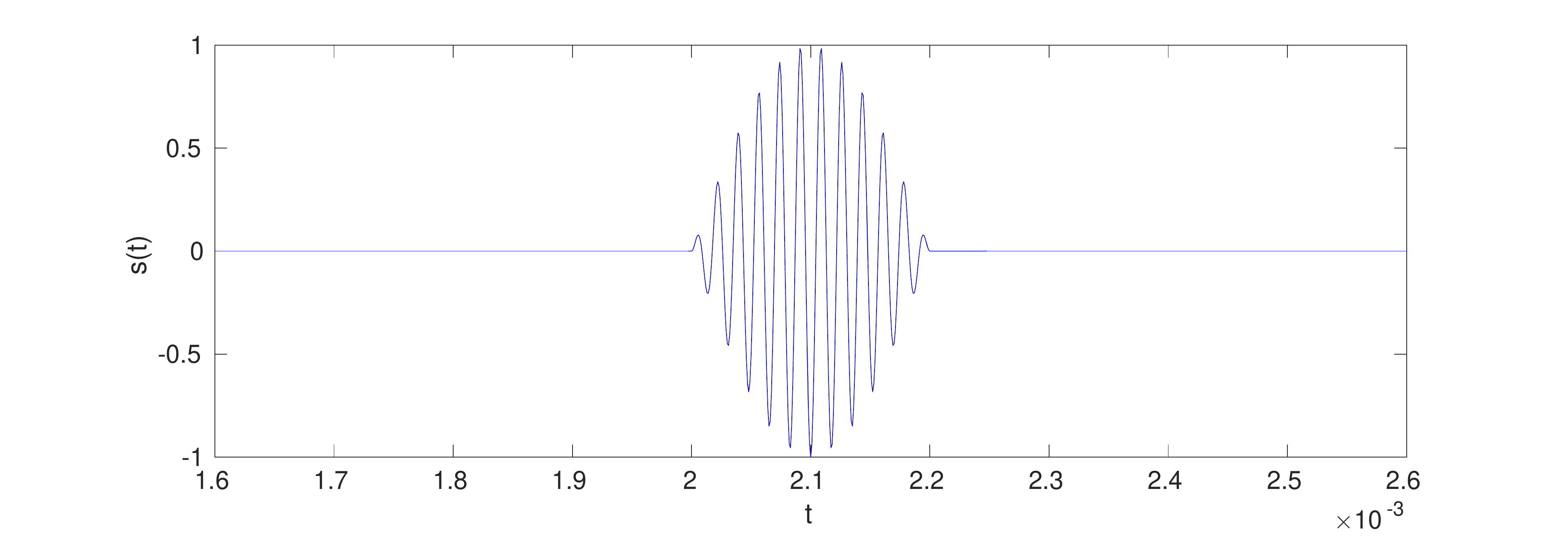}
    \caption{The function $s(t)$ given in~\eqref{eq:s(t)}.}
    \label{fig:s(t)}
\end{figure}

The experiment begins by sending a short pulse from the location $\br_0$. The signal $s(t)$ (shown in Figure~\ref{fig:s(t)}) is given by
\begin{equation} \label{eq:s(t)}
	s(t) := \sin(\omega_R t) \times \sin( 5000\pi(t - t_0 )) \times \1 \left( 0 < t - t_0 < 5000 \right), \quad \text{with} \quad t_0 = 2\times10^{-3} \ \text{s}.
\end{equation}

\begin{figure}[h!]
    \hspace{-2cm} \includegraphics[scale=0.5]{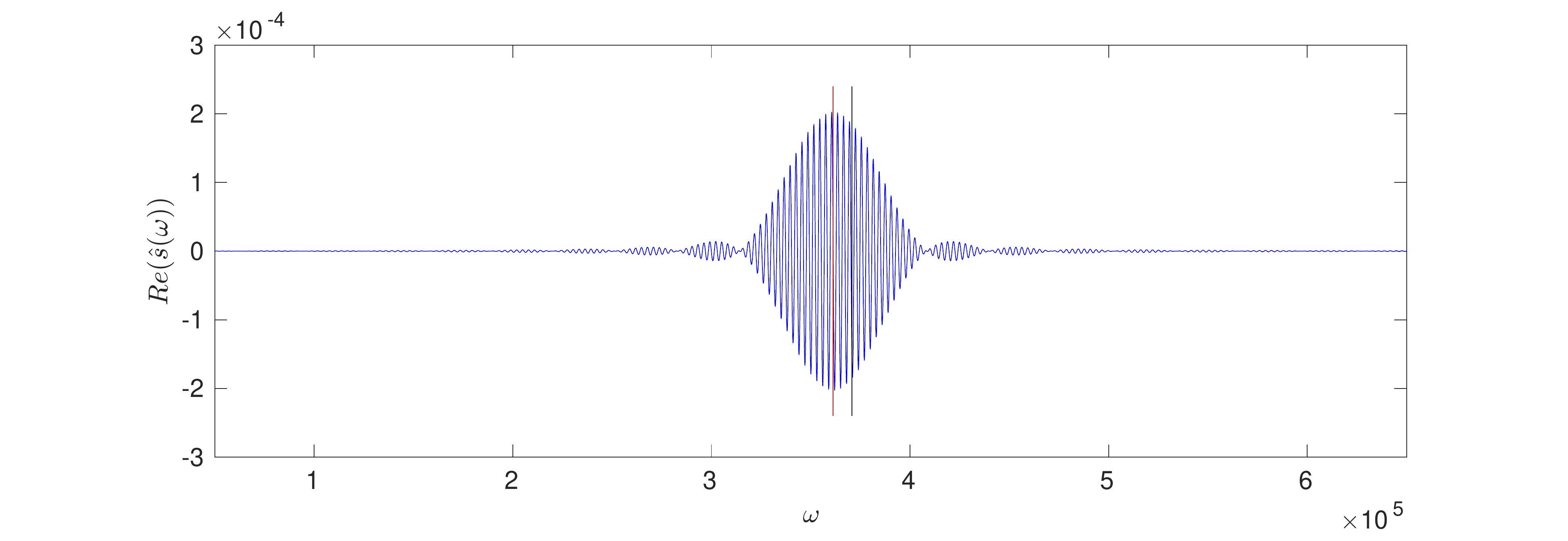}
    \caption{The function $\Re(\hat{s}(\omega))$ given in~\eqref{eq:ps}. The red line indicates the reference frequency $\omega_R$, and the black line is the Minnaert frequency $\omega_M$.}
    \label{fig:ps}
\end{figure}

To compute the forward problem, the signal is time-Fourier transformed with
\begin{equation} \label{eq:ps}
	s(t) = \int_{\R} \hat{s}(\omega) \re^{- \ri \omega t}.
\end{equation}
In practice, since we expect the recording to last much longer than the short impulse, the signal is recorded at each receiver for a time period of $[0, T]$, where $T = 50 \times 10^{-3} \ \text{s}$, with the numerical resolution of the signal being $\Delta t = 1\times10^{-6} \ \text{s}$. As a result, we are in possession of the Fourier transformed values for all $\omega \in [-\omega_{\rm max}/2, \omega_{\rm max}/2]$, with $\omega_{\rm max} \approx 3.1\times 10^{6} \ \text{Hz}$, and $\Delta \omega \approx 125 \ \text{Hz}$. The plot of $\hat{s}(\omega)$ is given in Figure~\ref{fig:ps}.

\begin{figure}[h!]
    \hspace{-2cm} \includegraphics[scale=0.5]{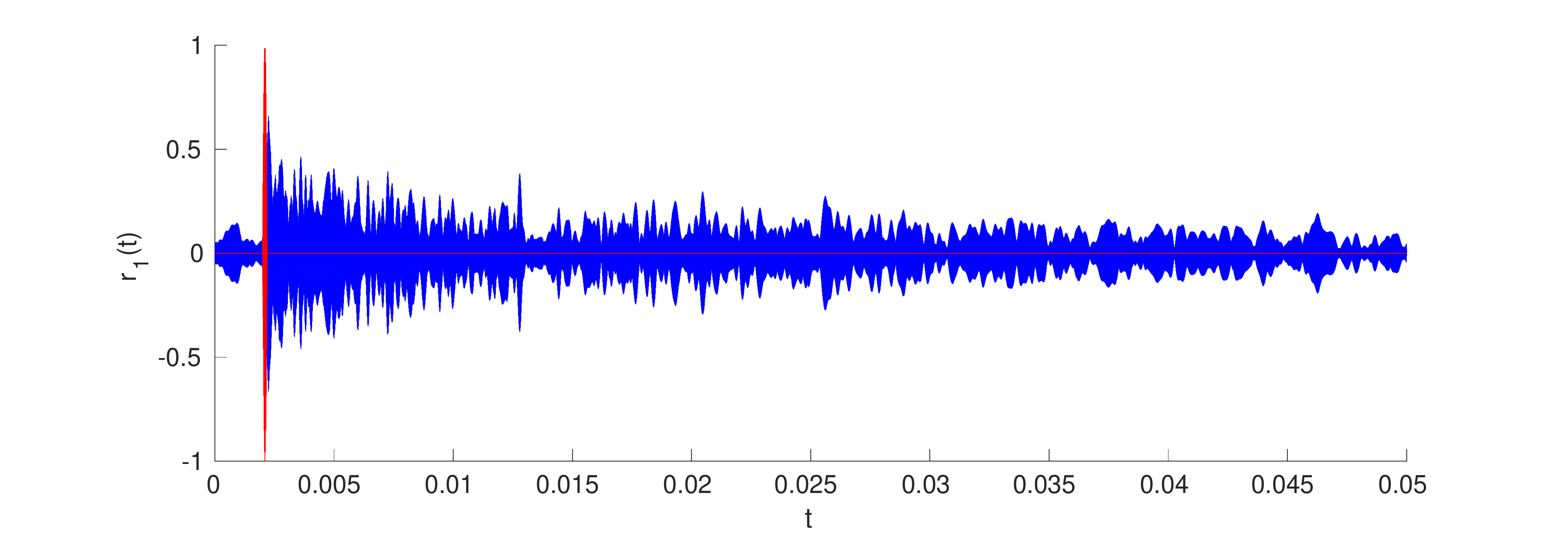}
    \caption{The recording at the first receiver $r_1(t)$ (blue). The red plot is $s(t)$.}
    \label{fig:r1(t)}
\end{figure}

Due to the fact that the signal is real-valued, and oscillates near $\omega_R$, only the frequencies near $\omega_R$ are considered (in practice this mean all frequencies in the range $\omega \in [0, 2 \omega_R]$). This necessitates solving approximately 5750 Helmholtz problems. We use the point scatterer model described in the previous section to solve these problems, and obtain the frequency domain representations of the recorded signals $\hat{r}_j(\omega)$ for $1 \le j \le 4$. Through application of the Fourier transform, we retrieve the recordings in the time domain $r_j(t)$ (see Figure~\ref{fig:r1(t)}). Note that while the initial signal is a short pulse, the recorded signal is much longer, which is due to the phenomenon of multiple echoing of waves among the bubbles. The non-null signal before $t_0$ is due to the Gibbs phenomenon.

\begin{figure}[h!]
    \hspace{-2cm} \includegraphics[scale=0.5]{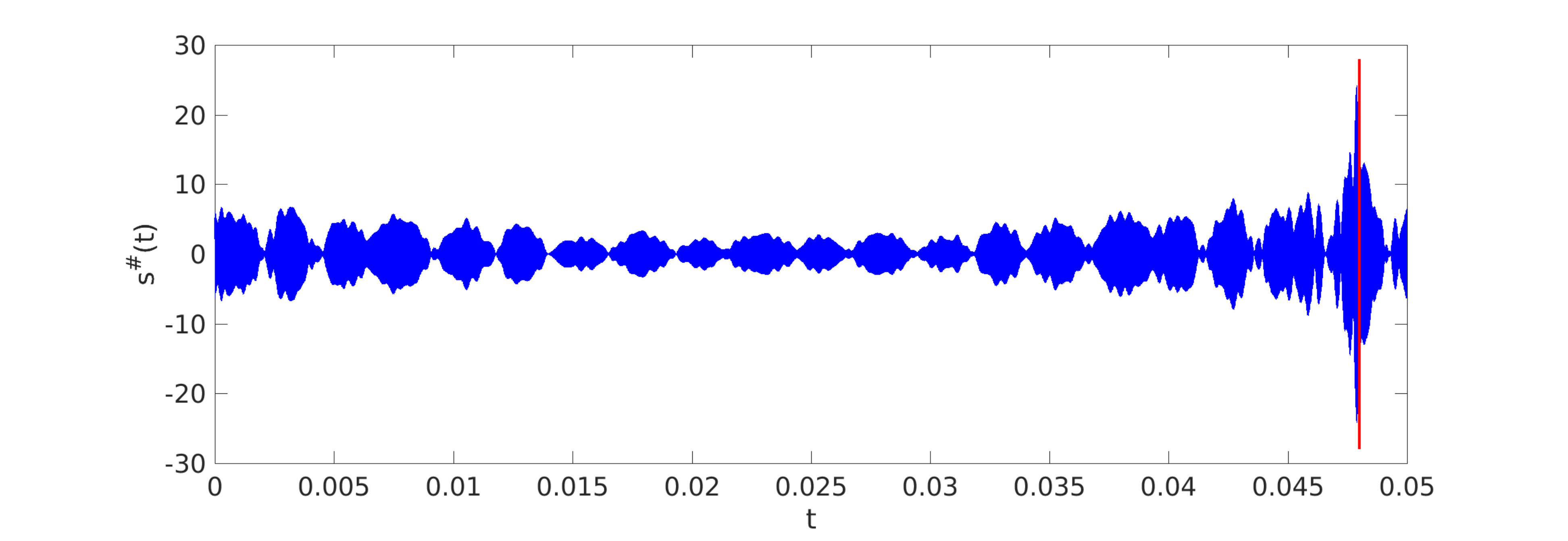}
    \caption{\scriptsize{The time-reversed signal $s^\sharp(t)$ at the origin. The red line indicates $T - t_0$, \textit{i.e.}, the time of the expected reversed pulse.}}
    \label{fig:Tr_t}
\end{figure}

\begin{figure}[h!]
	\renewcommand*\thesubfigure{\roman{subfigure}}
	\captionsetup[subfigure]{font=scriptsize}
    \begin{subfigure}{\linewidth}
	\hspace{-2cm} \includegraphics[scale=0.5]{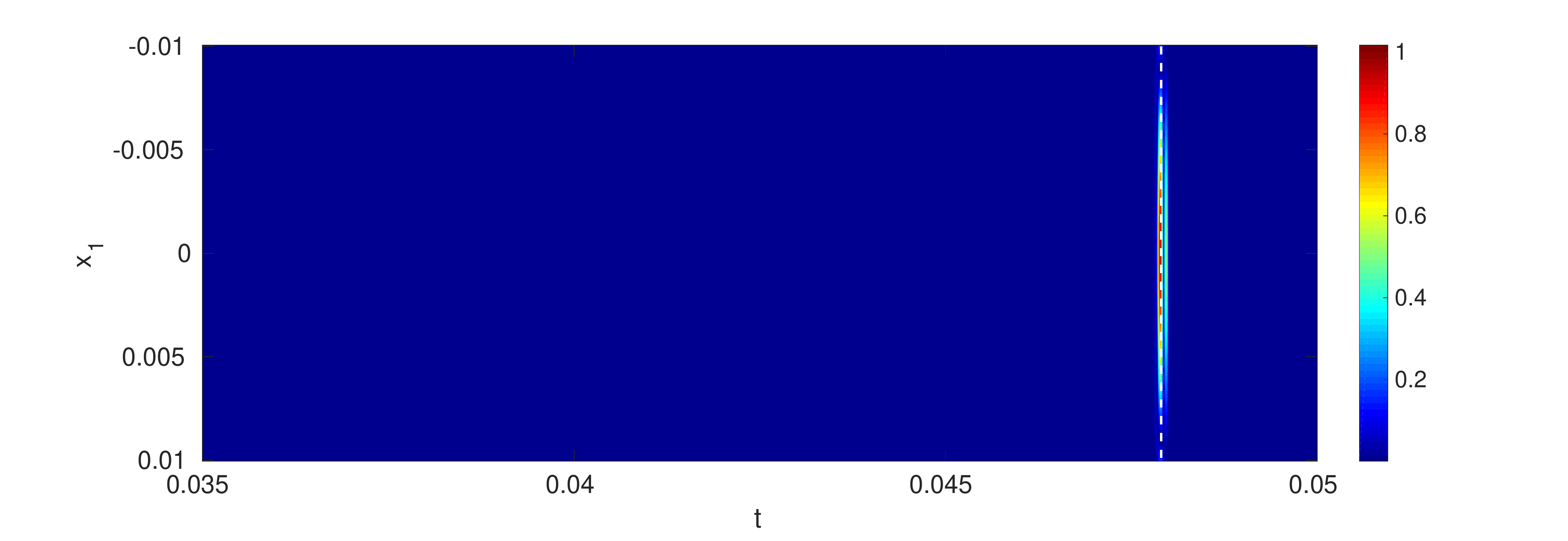}
        \caption{\scriptsize{The magnitude of the time-reversed pressure field $u(x,t)$ in the absence of bubbles.}}
     \end{subfigure}
	\begin{subfigure}{\linewidth}
    \hspace{-2cm} \includegraphics[scale=0.5]{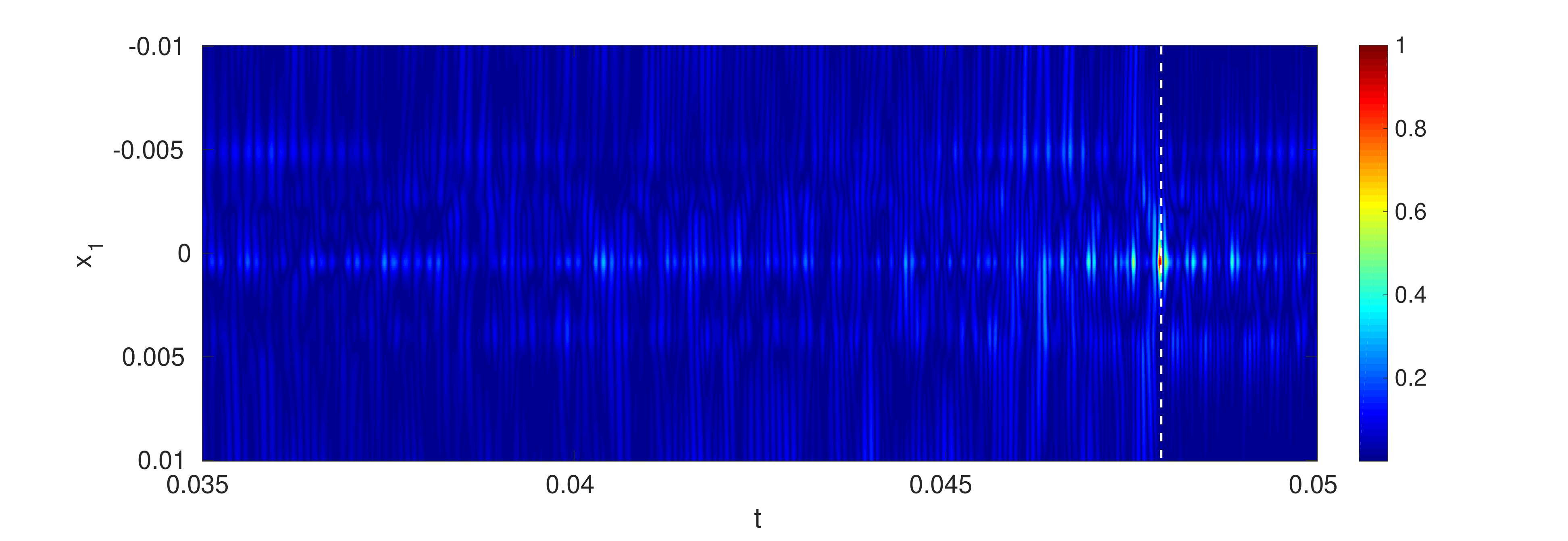}
        \caption{\scriptsize{The magnitude of the time-reversed pressure field $u(x,t)$ in the presence of bubbles.}}
     \end{subfigure}
        \caption{Magnitudes of normalized time-reversed pressure fields when bubbles are not present (i) and present (ii). The dashed white lines represent the expected arrival time $t=T-t_0 \approx 0.0480 \ \text{s}$ at the source.}
\label{fig:Tr_t_imagesc}
\end{figure}

We now simulate the time-reversal experiment. Each receiver now acts as a source, and emits the signal $r^\sharp_j(t) = r_j(T - t)$. Proceeding as in the forward problem, we time-Fourier transform the signal, and then solve the resulting Helmholtz problems. We record the time-reversed signal $s^\sharp(t)$ at the origin, and the result is shown in in Figure~\ref{fig:Tr_t}. In Figure~\ref{fig:Tr_t_imagesc} we plot the magnitudes of the normalized time-reversed pressure fields in the absence of bubbles, and in the presence of bubbles, respectively, for $x = (x_1,0,0)^\top$ where $x_1 \in [-L,L]$, and $t \in [0.035, 0.05]$. While in both plots it is clear that the time-reversed signal arrives at the origin at the expected time, greatly enhanced spatial localization (super-focusing in space) is readily apparent in the presence of bubbles.

%
%
%
\section{The imaginary part of the Green's function} \label{sec_imag_G}
It is known in imaging that the level of focusing achievable is linked to the imaginary part of the Green's function~\cite{ammari2013mathematical, Ammari2015_a}. In this section we concentrate on investigating the properties of the imaginary part of the Green's function when a large number of small acoustic bubbles are introduced to the region under consideration.

As in the previous section, we consider a source at the origin surrounded by a cube of uniformly distributed bubbles, in which the gas volume fraction is $\Phi = 2\times 10^{-4}$. We denote by $G_{\text{m}}(x, \omega)$ the Green's function in the presence of the bubbles at frequency $\omega$. Using the point interaction approximation, we have
\begin{equation} \label{eq:meanG}
	G_{\text{m}}(x, \omega) := G(x, \omega/c_w) + \sum_i 4 \pi f_s G(x - x_i, \omega/c_w) \sum_j M_{ij}^{-1} G(x_j, \omega/c_w),
\end{equation}
where $f_s$ is the scattering function defined in (\ref{eq:def:fs}),
$c_w$ is the wave speed in water and
\begin{equation}
M_{ij} =
\begin{cases}
1, & i = j, \\
4 \pi f_s G(x_i-x_j, \omega/c_w), & i \neq j.
\end{cases}
\end{equation}



\begin{figure}[ht!]
	\vspace{-0.0cm}  \hspace{-1.6cm}   \includegraphics[scale=0.46]{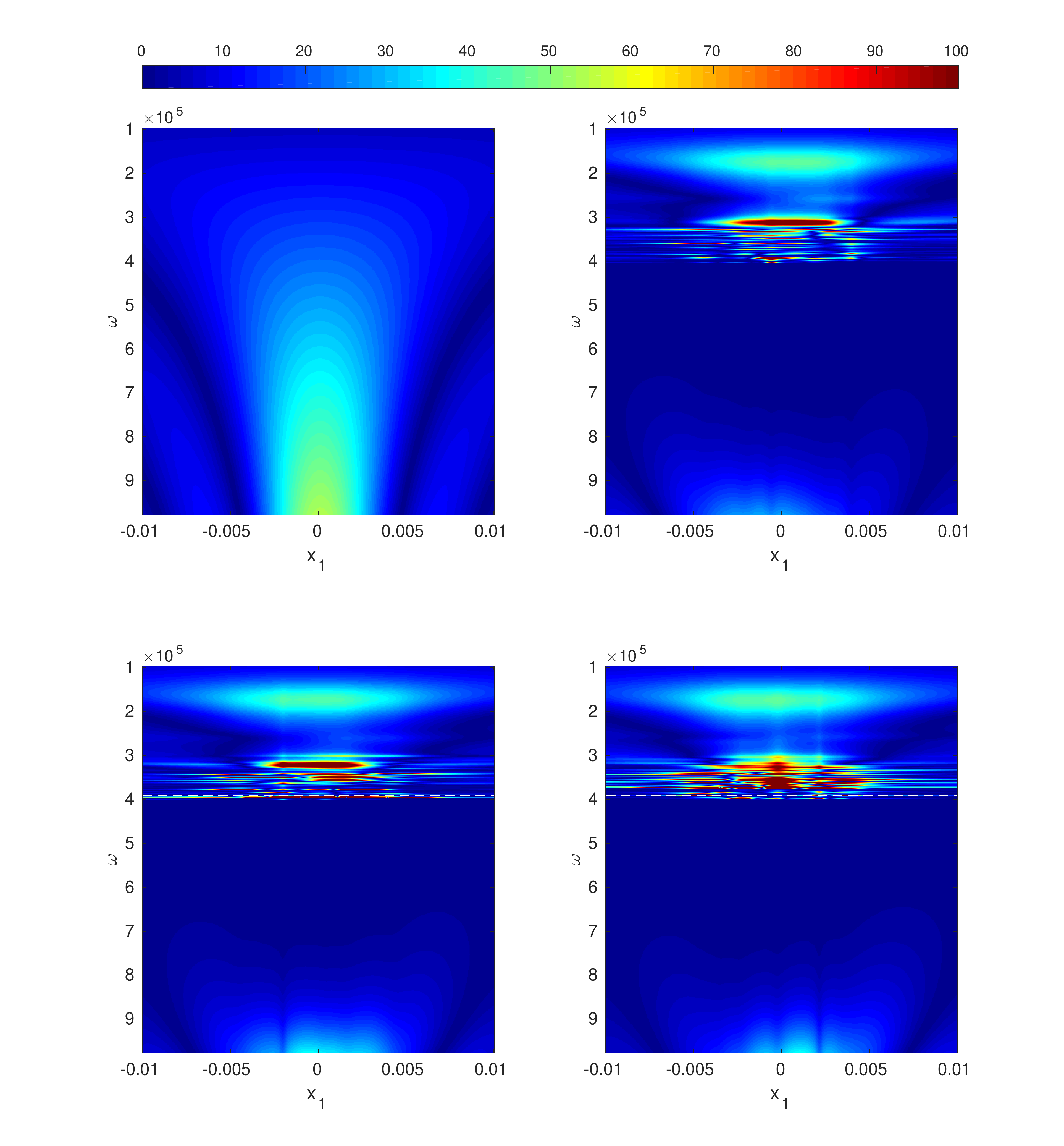}
    \caption{\scriptsize The magnitudes of the imaginary part of the free space Green's function $G(x,\omega)$ (top left) and the imaginary parts of the Green's functions for three different random configurations of bubbles $G_{\text{m}}(x,\omega)$ due to a source located at the origin. The dashed white line represents the Minnaert resonance frequency $\omega_M$.}
    \label{fig:imG}
\end{figure}

We set $\omega_{\rm min} = 15 \times 2 \pi \ \text{kHz}$ and $\omega_{\rm max} = 155 \times 2 \pi \ \text{kHz}$. In Figure~\ref{fig:imG}, we plot $|\Im(G(x, \omega))|$, and $|\Im(G_{\text{m}}(x, \omega))|$, for $x = (x_1,0,0)^\top$ where $x \in [-L,L]$ and $\omega \in [\omega_{\rm min}, \omega_{\rm max}]$ for three different random bubble configurations. We clearly observe that there are a lot of small oscillations with regard to the maps $\omega \mapsto \Im(G(x, \omega))$. This phenomenon was noticed and explained in the propagation of waves in random media, in the paraxial regime~\cite{Bal2002,Borcea2003}. We note that for each random configuration of bubbles, when the frequency is slightly less than the Minnaert resonance frequency (white dashed line), on average there is a very significant enhancement of the imaginary part of the effective Green's function, as compared to the free-space Green's function. Slightly above the Minnaert resonance frequency, on the other hand, the imaginary part of the Green's function vanishes almost completely.

\begin{figure}[!ht]
    \renewcommand*\thesubfigure{\roman{subfigure}}
	\captionsetup[subfigure]{font=scriptsize}
    \centering 
\begin{subfigure}{0.43\textwidth}
  \includegraphics[width=\linewidth]{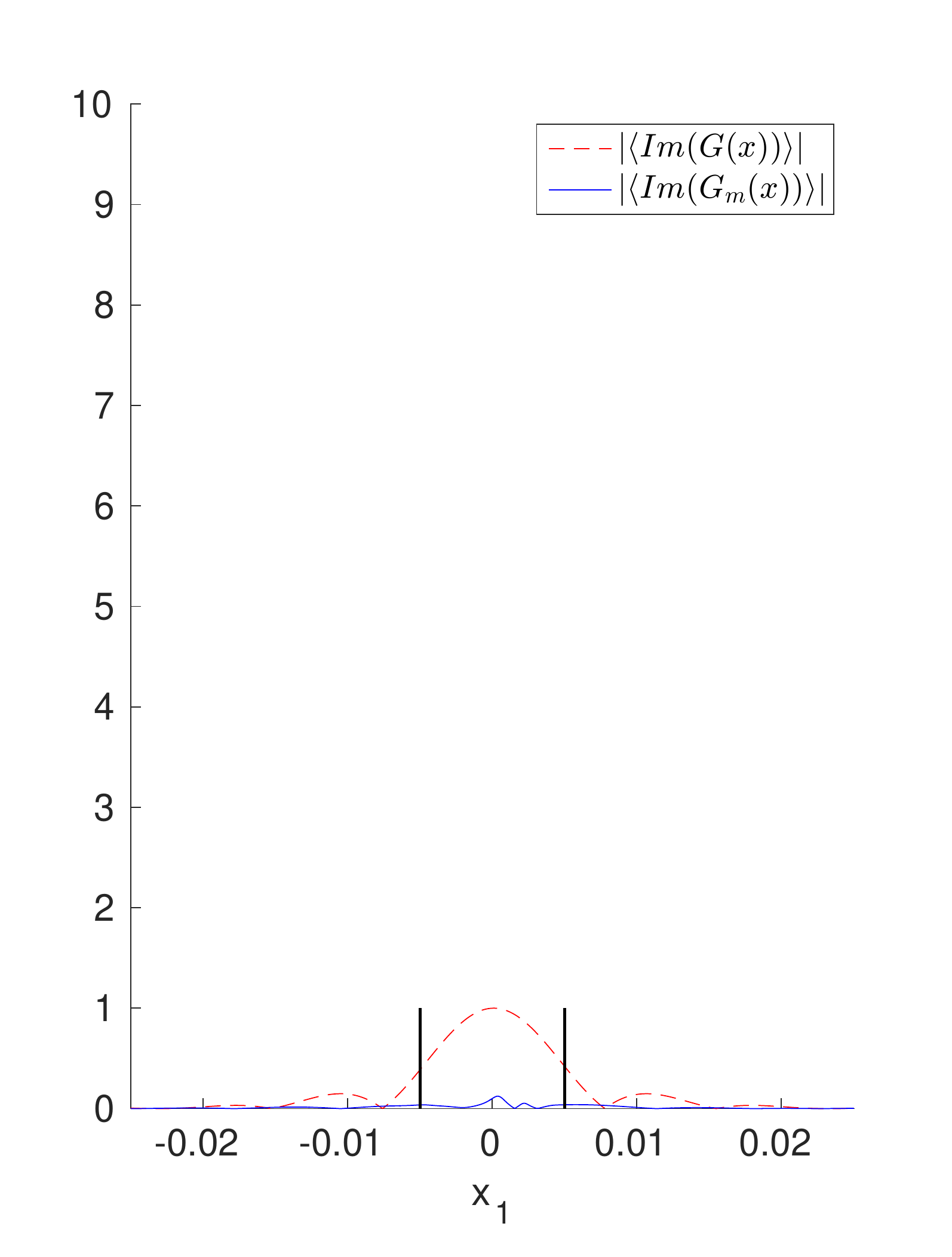}
  \caption{Averaging just above $\omega_M$.}
\end{subfigure}\hfil 
\begin{subfigure}{0.43\textwidth}
  \includegraphics[width=\linewidth]{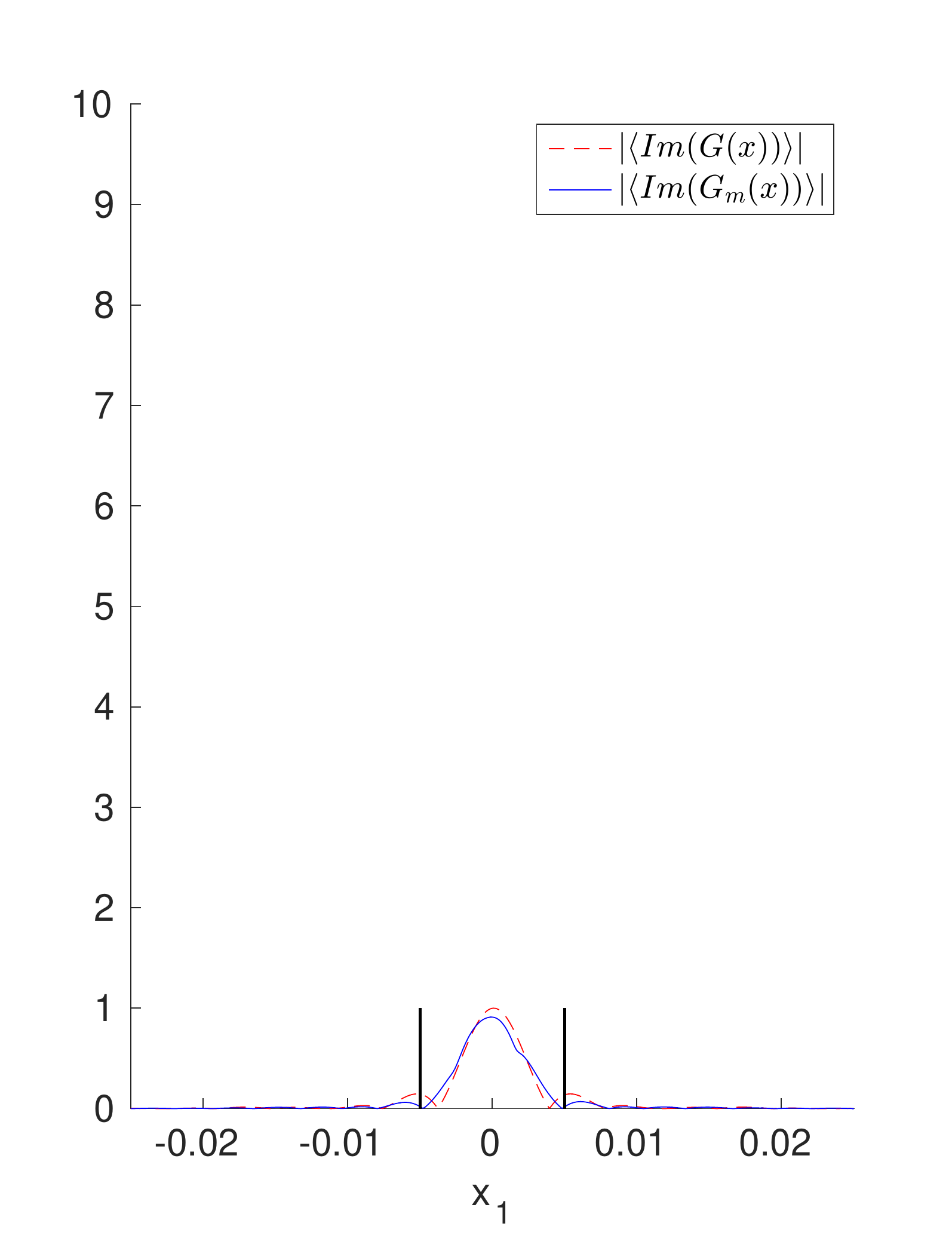}
  \caption{Averaging  well away from $\omega_M$.}
\end{subfigure}\hfill
\begin{subfigure}{0.43\textwidth}
  \includegraphics[width=\linewidth]{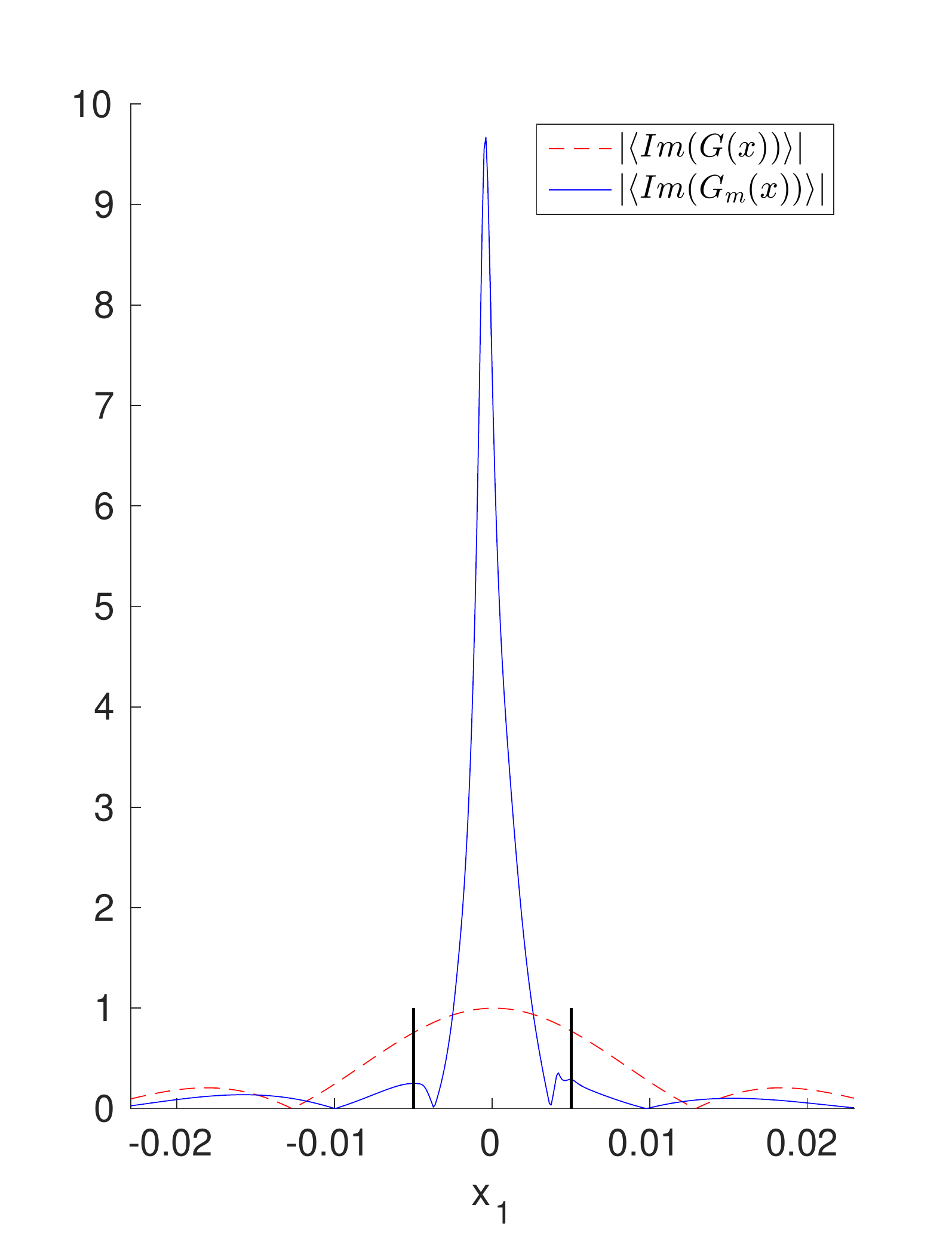}
  \caption{Averaging just below $\omega_M$.}
\end{subfigure}
        \caption{The normalized magnitudes of the averaged imaginary part of the free space Green's function $G(x)$ (red dashed line) and the averaged imaginary part of the Green's function in the presence of bubbles $G_{\text{m}}(x)$ (blue line), due to a source located at the origin. The averaging takes place over a range of frequencies just above the Minnaert resonance $\omega_M$ in (i), (above and) well away from  $\omega_M$ in (ii), and just below $\omega_M$ in (iii). The black lines show the limit of the box $[-L/2, L/2]$.}
\label{fig:avg_eff_imG}
\end{figure}

Finally, we consider the average of the imaginary part of the  Green's function over the frequency interval $[\omega_{\rm min}, \omega_{\rm max}]$:
\begin{equation} \label{eq:meanG2}
	\langle \text{Im} (G_{\text{m}}(x)) \rangle := \frac{1}{\omega_{\rm max} - \omega_{\rm min}} \int_{\omega_{\rm min}}^{\omega_{\rm max}} \text{Im}(G_{\text{m}}(x, \omega)) \rd \omega.
\end{equation}
In Figure~\ref{fig:avg_eff_imG} we compare the magnitude of this expression with the magnitude of the corresponding expression for the free-space Green's function. We consider three separate cases. In the first case the averaging is performed over a range of frequencies just above the Minnaert resonance frequency, \textit{i.e.}, $\omega_{\rm min} = 1.01 \  \omega_M$ and $\omega_{\rm max} = 2 \ \omega_M$. In this case the presence of the bubbles can be seen to have a strong effect on propagation and $\langle \text{Im} (G_{\text{m}}(x)) \rangle$ vanishes almost completely. Next we perform the averaging well away from and above the Minnaert resonance frequency, \textit{i.e.}, $\omega_{\rm min} = 2 \ \omega_M$ and $\omega_{\rm max} = 4 \ \omega_M$. It is clear that in this case the effect of the bubbles is greatly diminished and $\langle \text{Im} (G_{\text{m}}(x)) \rangle$ appears almost the same as it would in a bubble free system.

Finally, we perform the averaging over a range of frequencies just below the Minnaert resonance frequency, \textit{i.e.}, we set $\omega_{\rm min} = 0.8 \ \omega_M$ and $\omega_{\rm max} = 0.99 \ \omega_M$. Here, we observe a strong peak near the origin, which is indicative of the super-focusing phenomenon. The full width half maximum (FWHM) of the free space Green's function is $w \approx 0.0152$. In comparison, the FWHM of the Green's function in the presence of bubbles is $w_{m} \approx 0.0016$, indicating that the presence of a large number of small bubbles, which have been excited near their Minnaert resonance frequency, leads to a substantial increase in focusing power.


\subsection{Proof of Lemma~\ref{lem:SD_KD[1]}}
\label{sec:proof:SD_KD[1]}

This is a classical result, which is a special case of a more general result. It is usually proven using the properties of spherical Bessel functions (see for instance~\cite{Nedelec2001}). We include here a simple proof for completeness.

\medskip

~
 \textbf{Step 1}: Proof for $\cS_D^k$ with $x \in \partial D$.\\
Let us first prove the result for $\cS_D^k$. We have that for $x \in \partial D$,
\[
	\cS_D^k[\1_{\partial D}] (x) = - \dfrac{1}{4 \pi} \int_{\partial D} \dfrac{\re^{\ri k | x - y |}}{| x - y |} \rd \sigma (y).
\]
Since the problem is invariant under rotation, it is sufficient to compute the integral for $x = (0,0,R)$.
We use the spherical coordinates $$\partial D \ni y = (y_1, y_2, y_3)^T = (R \sin \phi \cos \theta , R \sin \phi \sin \theta, R \cos \phi)^T$$ with $0 \le \theta \le 2 \pi$ and $0 \le \phi \le \pi$. Then, for $y \in \partial D$, it holds that $$\left| x - y \right| = R \sqrt{ \sin^2 \phi + (1 - \cos \phi)^2} = {2}  R \sin \left( \frac{\phi}{2} \right),$$
so that
\[
	\cS_D^k[\1_{\partial D}] (x) = - \dfrac{2 \pi R^2}{4 \pi} \int_0^\pi \dfrac{\re^{2 \ri k R \sin \left( \frac{\phi}{2} \right)} }{2 R \sin \left( \frac{\phi}{2} \right)} \sin\left( \phi \right) \rd \phi.
\]
Together with the equality $\sin(\phi) = 2 \sin(\phi/2) \cos(\phi/2)$ and the change of variable $u = \sin(\phi/2)$, we end up with
\begin{equation} \label{eq:SD[1]}
	\cS_D^k[\1_{\partial D}] (x) = - R \int_0^1 \re^{\ri k {2} Ru} \rd u = - R  \left( \dfrac{\re^{\ri k {2} R} - 1}{2 \ri k R} \right) = - \dfrac{\re^{\ri k R}}{k} \sin \left( kR \right).
\end{equation}

\indent \textbf{Step 2}: Exact expression of $\cS_D^k[\1_{\partial D}]$ with $x \in \partial D$.\\
From the previous step, we can easily deduce an expression of $\cS_D^k[\1_{\partial D}]$ on $\R^3 \setminus D$. Indeed, the function $\cS_D^k[\1_{\partial D}]$ outside $D$ is the (unique) solution to the Dirichlet equation
\[
	\left\{ \begin{array}{l}
		(\Delta + k^2) u = 0, \\
		u \Big|_{\partial D} =  \cS_D^k[\1_{\partial D}] =  - \dfrac{\re^{\ri k R}}{k} \sin ( k R),
	\end{array} \right.
\]
together with the Sommerfeld radiation condition. The solution is therefore
\begin{equation} \label{eq:explicitS1}
	\cS_D^k[\1_{\partial D}](x) = - \dfrac{R \sin ( k R)}{k} \dfrac{ \re^{\ri k | x |}}{ | x |}.
\end{equation}

\indent \textbf{Step 3}: Proof for $\left( \cK_D^k \right)^*$.\\
For all $x \in \partial D$,
\[
	\left(\cK_D^k\right)^*[\1_{\partial D}](x) = \f{\p}{\p \nu_x} \cS_D^k[\1_{\partial D}] - \dfrac{1}{2}.
\]
Together with~\eqref{eq:explicitS1}, we obtain
\[
	\left(\cK_D^k\right)^*[\1_{\partial D}](x) = - \dfrac{R \sin(kR)}{k} \re^{\ri k R} \left( - \frac{1}{R^2} + \dfrac{\ri k}{R} \right) - \frac12
	= \dfrac{\sin(kR)}{kR} \re^{\ri k R} - \dfrac{\re^{2 \ri k R}}{2}.
\]

\newpage

\end{document}